\documentclass[10pt]{article}
\usepackage{geometry}                % See geometry.pdf to learn the layout options. There are lots.
\usepackage[parfill]{parskip}    % Activate to begin paragraphs with an empty line rather than an indent
\usepackage{graphicx}
\usepackage{amsmath, latexsym, amssymb}
\input xypic
\begin{document}
\title{The existence of  closed  3-forms  of $\tilde G_2$-type  on  7-manifolds}
\author{\bf  H\^ong-V\^an L\^e }
\date{}
\maketitle

\newcommand{\R}{{\mathbb R}}
\newcommand{\C}{{\mathbb C}}
\newcommand{\F}{{\mathbb F}}
\newcommand{\Z}{{\mathbb Z}}
\newcommand{\N}{{\mathbb N}}
\newcommand{\Q}{{\mathbb Q}}
\newcommand{\Hq}{{\mathbb H}}
\newcommand{\Aa}{{\mathcal A}}
\newcommand{\Bb}{{\mathcal B}}
\newcommand{\Cc}{{\mathcal C}}    %configuration space
\newcommand{\Dd}{{\mathcal D}}
\newcommand{\Ee}{{\mathcal E}}
\newcommand{\Ff}{{\mathcal F}}
\newcommand{\Gg}{{\mathcal G}}    %gauge transformations
\newcommand{\Hh}{{\mathcal H}}
\newcommand{\Kk}{{\mathcal K}}
\newcommand{\Ii}{{\mathcal I}}
\newcommand{\Jj}{{\mathcal J}}
\newcommand{\Ll}{{\mathcal L}}    %Loop space
\newcommand{\Mm}{{\mathcal M}}    %moduli space
\newcommand{\Nn}{{\mathcal N}}
\newcommand{\Oo}{{\mathcal O}}
\newcommand{\Pp}{{\mathcal P}}
\newcommand{\Qq}{{\mathcal Q}}
\newcommand{\Rr}{{\mathcal R}}
\newcommand{\Ss}{{\mathcal S}}
\newcommand{\Tt}{{\mathcal T}}
\newcommand{\Uu}{{\mathcal U}}
\newcommand{\Vv}{{\mathcal V}}
\newcommand{\Ww}{{\mathcal W}}
\newcommand{\Xx}{{\mathcal X}}
\newcommand{\Yy}{{\mathcal Y}}
\newcommand{\Zz}{{\mathcal Z}}
\newcommand{\zt}{{\tilde z}}
\newcommand{\xt}{{\tilde x}}
\newcommand{\Ht}{\widetilde{H}}
\newcommand{\ut}{{\tilde u}}
\newcommand{\Mt}{{\widetilde M}}
\newcommand{\Llt}{{\widetilde{\mathcal L}}}
\newcommand{\yt}{{\tilde y}}
\newcommand{\vt}{{\tilde v}}
\newcommand{\Ppt}{{\widetilde{\mathcal P}}}

\newcommand{\Remark}{{\it Remark}}
\newcommand{\Proof}{{\it Proof}}
\newcommand{\ad}{{\rm ad}}
\newcommand{\Om}{{\Omega}}
\newcommand{\om}{{\omega}}
\newcommand{\eps}{{\varepsilon}}
\newcommand{\Di}{{\rm Diff}}
\newcommand{\Pro}[1]{\noindent {\bf Proposition #1}}
\newcommand{\Thm}[1]{\noindent {\bf Theorem #1}}
\newcommand{\Lem}[1]{\noindent {\bf Lemma #1 }}
\newcommand{\An}[1]{\noindent {\bf Anmerkung #1}}
\newcommand{\Kor}[1]{\noindent {\bf Korollar #1}}
\newcommand{\Satz}[1]{\noindent {\bf Satz #1}}

\newcommand{\gl}{{\frak gl}}
\renewcommand{\o}{{\frak o}}
\newcommand{\so}{{\frak so}}
\renewcommand{\u}{{\frak u}}
\newcommand{\su}{{\frak su}}
\newcommand{\ssl}{{\frak sl}}
\newcommand{\ssp}{{\frak sp}}

\newcommand{\Cinf}{C^{\infty}}
\newcommand{\CS}{{\mathcal{CS}}}
\newcommand{\YM}{{\mathcal{YM}}}
\newcommand{\Jreg}{{\mathcal J}_{\rm reg}}
\newcommand{\Hreg}{{\mathcal H}_{\rm reg}}
\newcommand{\SP}{{\rm SP}}
\newcommand{\im}{{\rm im}}

\newcommand{\inner}[2]{\langle #1, #2\rangle}    %inner products
\newcommand{\Inner}[2]{#1\cdot#2}
\def\NABLA#1{{\mathop{\nabla\kern-.5ex\lower1ex\hbox{$#1$}}}}
\def\Nabla#1{\nabla\kern-.5ex{}_#1}

\newcommand{\half}{\scriptstyle\frac{1}{2}}
\newcommand{\p}{{\partial}}
\newcommand{\notsub}{\not\subset}
\newcommand{\iI}{{I}}               %unit interval [0,1]
\newcommand{\bI}{{\partial I}}      %boundary of same
\newcommand{\LRA}{\Longrightarrow}
\newcommand{\LLA}{\Longleftarrow}
\newcommand{\lra}{\longrightarrow}
\newcommand{\LLR}{\Longleftrightarrow}
\newcommand{\lla}{\longleftarrow}
\newcommand{\INTO}{\hookrightarrow}

\newcommand{\Sy}{\text{ Diff }_{\om}}
\newcommand{\Ex}{\text{Diff }_{ex}}
\newcommand{\jdef}[1]{{\bf #1}}
\newcommand{\QED}{\hfill$\Box$\medskip}

\newcommand{\UuU}{\Upsilon _{\delta}(H_0) \times \Uu _{\delta} (J_0)}
\newcommand{\bm}{\boldmath}

\medskip

\begin{abstract}
In this note we construct a first example of a closed 3-form of $\tilde G_2$-type
on $S^3\times S^4$. We prove that $S^3\times S^4$  does not admit  a homogeneous  
 3-form of  $\tilde G_2$-type. Thus our example is a first example of a closed 3-form of $\tilde G_2$-type on a compact 7-manifold
 which is not stably homogeneous.
\end{abstract}


\section {Introduction.}

Let $\Lambda ^k (V ^n)^*$ be the space of k-linear anti-symmetric forms
on a given linear space $V^n$. For each $\om\in \Lambda ^k (V^n)^*$ we denote by
$I_\om$ the linear map
$$ I_\om : V^n \to  \Lambda ^{k-1} (V^n)^*,\: x\mapsto (x\rfloor \om) \, : = \om ( x, \cdots ).$$

 A k-form $\om$ is called {\bf multi-symplectic}, if $I_\om$ is a monomorphism. Two k-forms
 are equivalent, if they are in the same orbit of the action of $Gl(V^n)$ on $\Lambda^k(V^n)^*$.

The classification  of multi-symplectic 3-forms on $\R^7$ has been done
by Bures and Vanzura [1].  There are together 8 types of these forms, among them
there are two  stable forms  of $G_2$-type  and $\tilde G_2 $-type. 
More precisely the orbits  of these stable  3-forms  under the action of $Gl(V^7)$
are open sets in $\Lambda ^3(V^7)^*$ and their  corresponding isotropy groups are the compact group $G_2$ and
its  dual non-compact group $\tilde G_2$.

There are  many  known results on  7-manifolds  admitting a (closed) 3-form  of $G_2$-type,  see e.g.   [2],[3], [6],[7], [9],  [10], [11]. Manifolds which admit a (closed) 3-form of  $\tilde G_2$-type
are less known. In particular, known examples of closed 3-forms of $\tilde G_2$-type  on compact 7-manifolds  up to now are  homogeneous examples
or  obtaining from those by  adding a small closed 3-form. We call such a closed 3-form of $\tilde G_2$-type stably homogeneous.  
  
In this note we construct a first example of a closed 3-form $\om^ 3$ of $\tilde G_2$-type on a 
manifold $X^7 = S^3 \times  S^4$ by identifying $X^7$ with a submanifold of the group
$SU(3)$ provided with the Cartan 3-form (Theorem 2.1.) In section 3 we prove that the codimension
of the action of the full automorphism group of $(X^7, \om^3)$ on $X^7$ is 1 (Theorem
3.1), thus $(X^7, \om^3)$ is not homogeneous.  Moreover we  prove that $X^7$ admits no homogeneous $\tilde G_2$-structure (Proposition 3.2).
So our example is a $\tilde G_2$-structure on $S^3 \times S^4$ with  maximal symmetrie.

This note also contains an appendix which gives a necessary and sufficient condition for
a closed 7-manifold to admit a $\tilde G_2$-structure. As a corollary we obtain many examples of
open 7-manifolds admitting a closed 3-form of $\tilde G_2$-type.

\medskip

\section {New example of a closed 3-form of $\tilde G_2$-type on $S^3\times S^4$.}

On each semi-simple Lie group $G$ there exists a natural bi-invariant 3-form $\phi^3$ which is defined
 at the Lie algebra $g= T_eG$ as follows
$$\phi^3  ( X, Y, Z) = <X, [Y, Z]>,$$
where $<,>$ denotes the Killing form on $g$. This 3-form $\phi^3$ is also called
the Cartan 3-form.

We claim that the Cartan 3-form $\om^3$ is multi-symplectic. To show the injectivity of the
linear map $I_{\phi^3}$ we notice that if $X \in  \ker I_{\phi^3}$ , then
$$< X, [Y,Z] >= 0 \text { for all }  Y,Z \in g.$$
But this condition is incompatible with the semi-simplicity of $g$.
\medskip

Let us consider  the group $G = SU(3)$. For each $1\le i\le j \le 3$
let $g_{ij}(g)$ be the complex function on $SU(3)$ induced  from the standard unitary
representation $\rho$ of $SU(3)$ on $\C^3$: $g_{ij} (g) : = <\rho(g)\circ e_i, \bar e_j>$. Here
$\{ e_1 = (1, 0, 0), \, e_2 = (0,1,0), \, e_3 = (0,0,1)\}$ is a unitary basis of
$\C^3$.  
Let us denote by $X^7$ the co-dimension 1 subset in $SU(3)$ which is defined by the equation
$Im (g_{11}(g) ) = 0$.

\medskip

{\bf 2.1. Theorem.} {\it  The subset $X^7$  is diffeomorphic to the manifold
$S^3\times S^4$. Moreover $X^7$ is provided
with a closed  3-form of $\tilde G_2$-type  which is the restriction of $\phi^3$ to $ X^7$.}

\medskip

{\it Proof.} Let $SU(2) $ be the subgroup in $SU(3)$ consisting of all $g\in SU(3)$ such that
$\rho (g) \circ e_1 = e_1$. We denote by $\pi$ the natural projection
$$\pi: SU(3) \to SU(3)/SU(2).$$
We identify  $SU(3)/SU(2)$ with the sphere $S^5 \subset \C^3$ using the standard representation
 $\rho$ of $SU(3)$ on $\C ^3$.  Namely we set
$$\tilde\rho (g\cdot SU(2)): =\rho( g )\circ e_1.$$

We denote by $\Pi$ the composition $\tilde\rho\circ \pi: SU(3) \to SU(3)/SU(2) \to S^5$.
Let $S^4 \subset S^5$ be the geodesic sphere which consists of  points $v\in S^5$ such that $ Im\, e^1(v) = 0$. Here $\{ e ^i, i =1,2, 3\}$ are the  complex 1-forms on $\C^3$
which are dual to $\{ e_i\}$.
The pre-image $\Pi ^{-1} (S^4)$  consists of all $g\in SU(3)$ such that
$$ Im \, e^1 (\rho(g) \circ  e_1)  = 0.$$
$$\LLR  Im\, (g_{11}) = 0.$$
So $X^7$ is $SU(2)$-fibration over $S^4$. But this fibration is the restriction
of the $SU(2)$-fibration $\Pi^{-1} ( D^5)$ over the half-sphere $D^5$ to the boundary
 $\p D^5 = S^4$. Hence $X^7$ is a trivial $SU(2)$-fibration.
This proves the first statement of Theorem 2.2.

\medskip

Let us denote by  $SO(2)^1$    the orthogonal group
of the real subspace $\R ^2 \subset \C^3$ such that $\R^2$ is the span of $e_1$ and
$e_2$ over $\R$. Clearly $SO(2)^1$ is a subgroup of $SU(3)$.
\medskip

We denote by $m_L (g)$ (resp. $m_R (g)$) the left multiplication (resp. the right multiplication) 
by  an element $g\in SU(3)$.

\medskip

{\bf 2.2. Lemma}.  {\it  $X^7$ is invariant under the action
of $m_L (SU(2))\cdot m_R (SU(2))$.  For each $v \in S^4$ there exist  an element $\alpha
\in SO(2)^1$ and an element $g \in SU(2)$ such that $\Pi (g \cdot \alpha) = v$. Consequently  for any point  $x\in X^7$ there are $g_1, g_2 \in SU(2)$ and $\alpha \in
SO(2) ^1$ such that}
$$x = g_1\cdot \alpha \cdot g_2, \leqno (2.2.1)$$

{\it Proof.} Using the identification
that $X^7 = \Pi^{-1} (S^4)$ we note that the orbits of $m_R (SU(2))$-action on $X^7$ are the fibers $\Pi ^{-1} (v)$. Now the first statement  follows by  straightforward calculations. Let $v = ( \cos \alpha, z_2, z_3)\in S^4$, where
$z_i \in \C$. We choose $\alpha \in SO(2)^1$ so that
$$\rho(\alpha) \circ e_1  = (\cos \alpha , \sin \alpha)\in \R^2. $$
Clearly $\alpha$ is defined by $v$ uniquely up to sign $\pm $. We set
$$ w : = ( \sin \alpha, 0) \in \C^2 = < e_2, e_3>_{\otimes \C}.$$
We notice that 
$$|z_2| ^2 + | z_3| ^2 = \sin ^2 \alpha. $$
Since $SU(2)$ acts transitively on the sphere $S^3 $ of radius $| \sin \alpha| $
in $ \C^2= <e_2, e_3> _{\otimes \C}$,  there
exists an element $g \in SU(2)$ such that $\rho(g) \circ w =  (z_2, z_3)$.
Clearly 
$$ \Pi (g \cdot \alpha) = v .$$
This proves the second statement.
The last statement of Lemma 2.2 follows from the second statement and  the fact that $X^7 = \Pi ^{-1} (S^4)$
\QED

\medskip

{\it  Continuation of the proof of the   second  part of Theorem 2.1.}  Let us  write here a canonical expression of  a 3-form $\om^3$ of $\tilde G_2$-type (see e.g. [2], [1])
$$\om^3=\theta_1 \wedge \theta_2 \wedge \theta_3 +\alpha_1 \wedge \theta _1 + \alpha_2 \wedge \theta_2 + \alpha_3\wedge \theta _3.\leqno(2.3)$$
Here $\alpha_i$ are 2-forms on $V^7$ which can be written as
$$\alpha_1 = y_1\wedge y_2 + y_3\wedge y_4, \: \alpha_2 = y_1 \wedge y_3 - y_2\wedge
y_4, \: \alpha_3 = y_1 \wedge y_4 + y_2 \wedge y_3$$
and $(\theta_1, \theta_2, \theta_3, y_1, y_2, y_3, y_4)$ is an  oriented basis of
$(V^7)^*$.

Using Lemma 2.2   we reduce the proof  of  the second statement of Theorem 2.1   to  verifying that the value of $\phi^3_{|X^7}$ at any point
$\alpha \in SO (2) ^1\subset X^7$ is a  3-form of $\tilde G_2$-type. 
First we shall compute that  value at point  $ e\in SO(2) ^1\subset  X^7$ and  then
 we shall compute that value   at any point  $\alpha \in SO(2) ^1$.
\medskip

\underline{ Step 1.} We shall use the Killing metric to identify the
Lie algebra $su(3)$  with its  co-algebra $g$.  In what follows we shall not distinguish co-vectors and vectors, poly-vector and exterior forms on $su(3)$. Clearly  we have
$$ T_e X^7 = \{ v \in su (3): Im \, g_{11} (v)  = 0\}.$$
Now we identify $gl(\C^3)$ with $\C^3 \otimes (\C^3) ^*$ and we denote
by $e_{ij}$ the element of $gl(\C^3)$ of the form $e_i \otimes (e_j)^*$.
Let $\delta _i$ be the  1-forms in $T_e X^7$ which are defined as follows:
$$ \delta_1 = {i\over \sqrt 2}(e_{22}- e_{33}), \delta _2 ={1\over \sqrt 2}( e_{23} - e_{32}), \delta_3 = {i\over \sqrt 2} (e_{23} +e_{32}).$$
\medskip

Furthermore, $\om_i$ are 2-forms on $T_eX^7$ which have the following expressions:
$$2 \om_1 = -(e_{12}-e_{21})\wedge i(e_{12}+ e_{21}) +  (e_{13}-e_{31})\wedge i(e_{13}+e_{31}),$$
$$2 \om_2 = -(e_{12}-e_{21}) \wedge  (e_{13}-e_{31})
- i(e_{12}+e_{21}) \wedge i(e_{13}+e_{31}),$$
$$2 \om_3 = -(e_{12}-e_{21})\wedge
i(e_{13}+ e_{31}) + i(e_{12}+e_{21})\wedge(e_{13}-e_{31}) . $$

{\bf 2.4. Lemma.} {\it The restriction of the Cartan 3-form to $X^7$ at the point $e$ is}
$$\phi^3_{|T_e X^7} = \sqrt 2 \delta _1 \wedge \delta_2 \wedge \delta_3 +{1\over  \sqrt2} \om _1 \wedge\delta _1 
+ {1\over \sqrt 2}\om _2\wedge \delta_2 + {1\over \sqrt 2}\om_3 \wedge \delta _3 ,\leqno(2.4.1)  $$

{\it Proof.} We set
$$f_1 := {e_{12}- e_{21}\over  \sqrt 2}
, \:  f_2 := {i(e_{12} + e_{21})\over 
\sqrt 2}$$
$$f_3 := {e_{13}-e_{31}\over \sqrt 2},
\: f4 := {i(e_{13} + e_{31})
\over \sqrt 2}$$
Then $(\delta_1, \delta_2, \delta_3, f_1, f_2, f_3, f_4$) form an ortho-normal basis of $T_eX^7$ w.r.t. to the restriction
of the Killing metric to $X^7$. Next we observe that $(\delta_1, \delta_2, \delta_3)$ form a $su(2)$-algebra.
Furthermore we have
$$ [\delta_2, \delta_3] = \sqrt 2\delta_1, [\delta_1, \delta_2] = \sqrt 2\delta_3, [\delta_1, \delta_3] = -
\sqrt 2\delta_2.\leqno(2.5.1)$$
$$ [f_1, f_2] = i(e_{11} - e_{22}), [f_1; f_3] = -{1\over
\sqrt 2}
\delta_2, [f_1, f_4] = -{1\over 
\sqrt 2\delta_3}.\leqno(2.5.2)$$
$$ [f_2, f_3] =
{1\over 
\sqrt 2}\delta_3,  [f_2, f_4] = -{1 \over 
\sqrt 2}\delta_2.\leqno(2.5.3)$$
$$ [f_3, f_4] = i(e_{11}- e_{33}).\leqno (2.5.4)$$
Using (2.5.1) - (2.5.4) we compute all the values $\phi^3(X_1,X_2,X_3)$ 
easily, where either $(X_1,X_2)$
are $(\delta_i, \delta_j)$ or $(X_1,X_2)$ are $(f_i, f_j)$ for some $(i, j)$ and $X_3$ is one of the basic vectors $(\delta_i, f_i)$.
In this way we get the equality (2.4.1). \QED

Now compare (2.4.1) with (2.3) we observe that these two 3-forms are $Gl(R^7)$ equivalent
(e.g. by rescaling $\delta_i$ with factor (1/2)). This proves that $\phi^3
_{|T_eX^7}$ is a 3-form of $\tilde G_2$-type.
 This completes the first step.

\medskip

\underline{ Step 2.} Using step 1 it suffices to show that 
$$D\,m_L(\alpha^{-1})( T_{\alpha} X^7) =  T_e X^7\leqno (2.6)$$ 
for any $ \alpha\in SO(2) ^1 \subset X^7$, $\alpha \not = e$.

Since   $X^7\supset \alpha\cdot SU(2)$, we have
$$su (2)\subset D\, m_L(\alpha^{-1})( T_{\alpha} X^7) .\leqno (2.7)$$
Denote by $SO(3)$ the standard orthogonal group of $\R^3 \subset \C^3$.
Since  $ \alpha \in SO(3) \subset X^7$, we have $D\, m_L (\alpha ^{-1} )( T_\alpha SO (3))
\subset D\, m_L (\alpha ^{-1} ) (T_\alpha X^7)$. In particular we have
$$ <(e_{12}-e_{21}), (e_{13}-e_{31})>_{\otimes \R} \subset D\, m_L(\alpha^{-1})( T_{\alpha} X^7).\leqno (2.8)$$

Since $SU(2)\cdot \alpha \subset X^7$, we have
$$ Ad (\alpha ^{-1}) su (2)  \subset D\, m_L(\alpha^{-1})( T_{\alpha} X^7).\leqno (2.9)$$

Using the formula 
$$Ad (\alpha^{-1}) = \exp ( -ad (t\cdot{e_{12}-e_{21}\over \sqrt 2})),\: t\not \equiv 0$$
we get immediately from (2.7), (2.8), (2.9) the following inclusion
$$ < i (e_{12} +e_{21}), i (e_{13} +e_{31})> _{\otimes \R} \subset D\, m_L (\alpha ^{-1}) (T_\alpha X^7))$$
which together with (2.7), (2.8) imply the desired  equality (2.6).

This completes the proof of Theorem 2.1.\QED

\medskip

\section{ The orbits of the action of the automorphism group of
$(X^7, \om^3)$ on $X^7$.}

Denote by $Aut(X^7, \om^3)$ the full automorphism group of the manifold $X^7$ equipped by the
3-form $\om^3$ constructed in Theorem 2.1. In this section we shall prove the following
\medskip

{\bf 3.1. Theorem.} {\it  The co-dimension of the $Aut(X^7, \om^3$)-action on $X^7$ is equal to 1.}

\medskip

Lemma 2.2 implies that the group $SU(2)\times SU(2)$ is a subgroup of the automorphism group
$Aut(X^7,\om^3)$. Further taking into account (2.2.1) we note that the dimension of a generic
orbit of the $Aut(X^7, \om^3)$-action on $X^7$ is at least 6. Hence Theorem 3.1 is a corollary of the
following Proposition, which implies that our manifold $(X^7, \om^3)$ is not stably homogeneous.

\medskip

{\bf 3.2. Proposition.}  {\it Let $G$ be a Lie group which acts transitively on $X^7$. Then there does
not exist a $G$-invariant 3-form of $\tilde G_2$-type on $X^7$.}

\medskip

{\it Proof of Proposition 3.2}. Since $S^3 \times S^4$ is connected, we can assume that $G$ is connected,
because the identity component of $G$ acts also transitively on $S^3 \times S^4$. Denote by $\bar G$
the simply connected covering of $G$. Then $\bar G$ acts transitively and almost effectively on
$X^7$.
We shall prove the following
\medskip

{\bf 3.3. Proposition.} {\it  Let $X^n$ be a compact connected space with $\pi_1(X^n) = 0 = \pi_2(X^n)$.
Suppose that $\bar G$ is a connected and simply connected group which acts transitively $X^n$. Let
$\bar G_{s,u}$ be a maximal compact group of $\bar G$. Then $\bar G_{s,u}$ acts transitively on $X^n$}.
\medskip

{\it Proof.} We denote by $\bar H$ the isotropy of the action of $\bar G$ on $X^n$. Since $X^n$ is simply
connected, the subgroup $\bar H$ is connected.
From the homotopy exact sequence
$$\pi_2(X^n) = 0 \to  \pi_1(\bar H ) \to  \pi_1(\bar G) = 0$$
we obtain that $\bar H $is also simply connected.
Using the Levi-Maltsev decomposition theorem we write
$$\bar G
= \bar G_s \times V^ k, \: \bar H = \bar H_s \times V ^r,$$
where $\bar G_s$ is a semisimple Lie subgroup of $\bar G$, $V^k$ is a solvable normal subgroup of $\bar G$ and
$\bar H_s$ is semisimple Lie subgroup of $\bar H$ and $V^r$ is a solvable normal subgroup of $\bar H$. Moreover
we can assume that $\bar H_s$ is a subgroup of $\bar G_s$ after applying an inner automorphism group
of $\bar G$.
Using the fbration
$$V^ r \to  \bar G/\bar H \to  \bar G/\bar H_s$$
we conclude that the quotient space $\bar G/\bar H_s$ is homotopic to $\bar G/\bar H$. Further using the fibration
$V^ k \to  \bar G/\bar H_s \to\bar  G_s/\bar H_s$ we conclude that $\bar G_s/\bar H_s$ has the same homotopy type of $X^7$.

Denote by $h_s$ the Lie algebra of $\bar H_s$ and by $h_{s,u}$ the maximal compact Lie sub-algebra of
the semi-simple Lie algebra $h_s$. The Lie subalgebra $h_{s,u}$ is the Lie algebra of a maximal
compact Lie group $\bar H_{s,u}$ of $\bar H_s$. We can also assume that $h_{s,u}$ is a subalgebra of a maximal
compact Lie subalgebra $g_{s,u}$ of $g_s$. Denote by $\bar G_{s,u}$ the maximal compact subgroup of $\bar G_s$
whose Lie algebra is $g_{s,u}$. Using the Iwasawa decomposition we conclude that the quotient
space $\bar G_{s,u}/\bar H_{s,u}$ has the same homotopy type of $X^n$.
Since $\bar G_{s,u}$ and $\bar H_{s,u}$ are compact and using
$$H_{i+1}(X^n,\Z) = 0, \text{  if } i \ge  n, \: H_n(X^n,\Z) = \Z,$$
we obtain that $\dim  \bar G_{s,u} = \dim \bar H_{s,u} + n$. Thus to show that the action of 
$\bar G_{s,u}$ on $X^n$ is
transitively, it suffces to prove
$$ g_{s,u }\cap h = h_{s,u},\leqno(3.4)$$
where $h$ is the Lie algebra of $\bar H$. We have the Iwasawa and Levi-Maltsev decomposition
$$h = h_{s,u} + h_f + r,$$
where $r$ is the radical of h which is the Lie algebra of $V^r$. Suppose that $v \in ((g_{s,u}\cap h)\setminus h_{s,u})$.
Since $h_{s,u} \in g_{s,u}$ we can assume that $v \in (h_f + r)$ after adding some vector $w$ in $h_{s,u}$, if
necessary. Let $v = v_f + v_r$ be the decomposition of $v$ into the components in $h_f$ and $r$.
Assume that $v_f \not= 0$. Then the the closure of the 1-parameter subgroup $\exp tv$ in $\bar H$ is noncompact,
since its projection to the quotient group $(\bar H /V^r) = \bar H_s$ is non-compact. Since $\bar H$
is a closed subgroup of $\bar G$, we conclude that the closure of $\exp tv$ in $\bar G$ is non-compact. On
the other hand, the subgroup $\exp tv$ lies in the compact subgroup $\bar G_{s,u}$. Thus we arrive at a
contradiction. Hence $v_f = 0$, so $v = v_r$. If $v_r \not= 0$, then the closure of the subgroup $\exp tv$
in $\bar H$ lies in $V^r$ and it is also non-compact. Hence $v_r$ also vanishes. This proves (3.4) and
completes the proof of Proposition 3.3.\QED
\medskip

{\it Continuation of the proof of Proposition 3.2.}  Using Proposition 3.2 it suffces to consider
only the transitive action of the compact, connected and simply connected group $\bar G_{s,u}$ on
$X^7 = S^3 \times S^4$. Next we can assume that this action is almost effective, since otherwise
the quotient $\bar G_{s,u}/N$ acts effectively and transitively on $X^7$, where $N$ is the kernel of the
action. The group $\bar G_{s,u}/N$ is compact, connected, and semi-simple. Its compact universal
covering acts almost effectively on $X^7$.

Next we observe that the isotropy group of the action of $\bar G_{s,u}$ on $X^7$ is $\bar H_{s,u}$, since the
isotropy group of this action must be connected. From now on we shall assume that $\bar G_{s,u}=\bar H_{s,u}$ admits a $\bar G_{s,u}$-invariant 3-form $\om^3$ of $\tilde G_2$-type
on $X^7$.

Let $e = e \cdot \{ H \}$ be a reference point on $X^7$. Denote by $\rho$ the isotropy representation of $h_{s,u}$
on $T_e( \bar G/\bar H ) = \R^7$. Since $h_{s,u}$ is semi-simple, the representation $\rho$ is faithful, i.e. the kernel
of $\rho$ is empty.
 Clearly $\rho(h_{s,u})$ must be a sub-algebra of the Lie algebra $\tilde g_2$ of $\tilde G_2$, in particular the
rank of $h_{s,u}$ is at most 2.
\medskip

{\bf 3.5. Lemma.} {\it  We have $h_{s,u} = su(2) \oplus su(2)$.}

\medskip

{\it Proof.} Since $\rho$ is faithful, we identify $h_{s,u}$ with its image in $\tilde g_2$. Clearly the algebra $h_{s,u}$ is
also a subalgebra of a maximal compact Lie subalgebra in $\tilde g_2$, so $h_{s,u}$ is either $su(2)$ or
$su(2) \times su(2)$.

Suppose that $h_{s,u} = su(2)$. Then $\dim \bar G_{s,u} = 10$ and therefore $\bar G_{s,u}$ must be a product
of classical Lie groups. In particular we know that $\pi_4(\bar G_{s,u})$ is a finite group (see e.g. [8,
\S 25.4].)

We also have $\pi_3(\bar H_{s,u}) = \Z$. Further we observe that the inclusion $\pi_3(\bar H_{s,u}) \to  \pi_3(\bar G_{s,u})$
 is
injective, since the subgroup $\bar H_{s,u}$ realizes a non-trivial element in $H_3(\bar G_{s,u}, \R)$.
Now let us consider the homotopy exact sequence
$$ \pi_4(\bar G) \to \pi_4(X^7) = \Z_2 \oplus  \Z \to  0 = \ ker(\pi_3(\bar H ) \to \pi_3(\bar G)). \leqno (3.6)$$
The exact sequence (3.6) implies that $\pi_4(\bar G)$ contains a subgroup $\bar Z$ which is impossible.
Hence $h_{s,u}$ cannot be $su(2)$.\QED

\medskip

{\bf 3.7. Lemma.} {\it  We have $g_{s,u} = su(2) \oplus so(5)$.}

\medskip

{\it Proof.} Since $\dim g_{s,u} = 6+7 = 13$, it is easy to see that $g_{s,u}$ must be the product of classical
compact Lie groups. Furthermore we notice that, since the dimension of $\bar G_{s,u}/\bar H_{s,u}$ is odd,
the rank of $\bar G_{s,u}$ must be strictly greater than the rank of $\bar H_{s,u}$. Let $g^1
_{s,u}$ be a simple
component of $g_{s,u}$. Then the rank of $g^1
_{s,u}$ is less than or equal to 2, otherwise the dimension
of $g_{s,u}$ is greater than or equal to 15. Thus $g_{s,u}$ must be a sum of simple components of
rank 1 or 2. Next, by dimension reason, $g_{s,u}$ cannot be a sum of only components of rank
1 and it cannot contain more than one simple component of rank 2. Looking at the table of
simple Lie algebras we arrive at  Lemma 3.7. \QED
\medskip

{\bf 3.8. Lemma.} {\it  The subalgebra $h_{s,u}$ lies in the component $so(5)$ of $g_{s,u}$.}

\medskip

{\it Proof.} We notice that $\bar G_{s,u}$ has the same homotopy type as $\bar G$ and hence it is simply
connected. Thus $\bar G_{s,u}$ is the product $SU(2)\times Spin(5)$. Analogously $\bar H_{s,u}$ must be 
$SU(2)\times
SU(2)$. Denote by $e_1$ the composition of the embedding $e : h_{s,u} \to  g_{s,u}$ with the projection
$p : g_{s,u} \to su(2)$. To prove Lemma 3.8 it suffces to show that the kernel of $e_1$ is equal to
$h_{s,u}$.

Suppose that the kernel of $e_1$ is not equal to $h_{s,u}$. Then this kernel must be one of the
components $su(2)$ of $h_{s,u}$ since $e_1$ can not be injective. Denote this kernel by $su(2)^2$ and
let $su(2)^1$ be the other component of $h_{s,u}$. The isomorphism $e_1 : su(2)^1 \to su(2)$ lifts to
an isomorphism denoted by $\tilde e_1$ from the corresponding component $SU(2)^1$ of $\bar H_{s,u}$ to the
component $SU(2)$ of $\bar G_{s,u}$. Denote by $\tilde e_2$ the homomorphism from $SU(2)^1 \times SU(2)^2$ to
$Spin(5)$. Since $\tilde e_1(SU(2)^2) = Id$, the restriction of $\tilde e_2$ to $SU(2)^2$ is an embedding. We shall
construct a map
$$I : \bar G_{s,u}/\bar H_{s,u} \to  Spin(5)/\tilde e_2(SU(2)^2)$$
and show that  $I$ is a homeomorphism. For each point $(a \cdot  b) \cdot \{ \bar H_{s,u}\} \in \bar G_{s,u}/\bar H_{s,u},$
where $a \in  SU(2)$ and $b \in  Spin(5)$ we set
$$I((a, b) \cdot  \{ \bar H_{s,u}\} ) = (b \cdot [\tilde e_2(\tilde e_1^{-1}(a))]^{-1}) \cdot  \{\tilde e_2(SU(2)^2)\} \in  Spin(5)/\tilde e_2(SU(2)^2).\leqno(3.9)$$
The map (3.9) is well-defined, since for any $(\theta_1,\theta_2) \in SU(2)^1 \times SU(2)^2$ we have
$$I(a \cdot  \tilde e_1(\theta_1), b \cdot  \tilde e_2(\theta_1, \theta_2))) \cdot \{ \bar H_{s,u}\}) = (b \cdot \tilde e_2(\theta_1, \theta_2) \cdot [\tilde e_2(\tilde e_1^{-1}(a \cdot \tilde e_1(\theta_1))]^{-1}) \cdot  \{\tilde e_2(SU(2)^2)\}$$
$$= (b \cdot  \tilde e_2(\theta_1,  \theta_2) \cdot  [\tilde e_2(\tilde e_1^{-1}(a) \cdot  \theta_1)]^{-1}) \cdot \{\tilde e_2(SU(2)^2)\}$$
$$= (b \cdot \tilde e_2(\theta_1, \theta_2) \cdot [\tilde e_2(\theta_1)]^{-1} \cdot [\tilde e_2(\tilde e_1^{-1}(a))]^{-1}) \cdot  \{\tilde e_2(SU(2)^2)\}$$
$$= (b \cdot  \tilde e_2(\theta_2) \cdot [\tilde e_2(\tilde e_1^{-1} (a))]^{-1}) \cdot \{\tilde e_2(SU(2)^2)\}$$
$$= (b \cdot [\tilde e_2(\tilde e_1^{-1} (a))]^{-1} \cdot  \tilde e_2(\theta_2)) \cdot \{\tilde e_2(SU(2)^2)\} = I((a, b) \cdot \{ \bar H_{s,u}\}).$$
Substituting $a = 1$ into (3.9) we obtain that $I$ is surjective. Now suppose that $I((a; b) \cdot 
\{\bar  H_{s,u}\}) = I((a', b') \cdot  \{ H_{s,u}\})$. Then
$$b' \cdot  [\tilde e_2(\tilde e_1^{-1}(a'))]^{-1} = b \cdot [\tilde e_2(\tilde e_1^{-1} (a))]^{-1}\cdot \tilde e_2(\theta_2) = b \cdot \tilde e_2(\theta_2) \cdot [\tilde e_2(\tilde e_1^{-1}(a))]^{-1},$$
for some $\theta_2 \in SU(2)^2$. Hence
$$ b' = b \cdot  \tilde e_2(\theta_2) \cdot  [\tilde e_2(\tilde e_1^{-1}(a))]^{-1} \cdot [\tilde e_2(\tilde e_1^{-1}(a'))] = b \cdot \tilde e_2(\theta_2) \cdot \tilde e_2[(\tilde e_1^{-1} (a^{-1})) \cdot  \tilde e_1^{-1}(a'))].\leqno(3.10)$$
Set $\theta_1 = \tilde e^{-1}(a^{-1} \cdot  a')$. Then  $a' = a \cdot  \tilde e_1(\theta_1)$,  and we get from (3.10)
$$b' = b \cdot  \tilde e_2(\theta_2) \cdot  \tilde e_2(\theta_1) = b \cdot  \tilde e_2(\theta_1) \cdot  \tilde e_2(\theta_2).$$
Hence $I$ must be injective. Thus we have proved that $I$ is a homeomorphism. Hence
$\pi_4(Spin(5)/\tilde e_2(SU(2)^2) = \Z$. Now let us consider the homotopy exact sequence
$$\pi_4(Spin(5))\stackrel{ j}{\to}\pi_4(Spin(5)/\tilde e_2(SU(2)^2)) \to \pi_3(\tilde e_2 (SU(2)^2) \stackrel{i}{\to} \pi_3(Spin(5)).$$
Clearly  $i$ is injective. Hence $j$ must be surjective. But $\pi_4(Spin(5))$ is a finite group.
Therefore $\pi_4(Spin(5)/\tilde e_2(SU(2)^2)$ is also a finite group. Thus we arrive at a contradiction, since
$\pi_4(X^7)$ contains an infinite subgroup. This proves Lemma 3.8. \QED
\medskip

{\it Completion of the proof of Proposition 3.2.} From Lemma 3.8 it follows that the restriction
of the isotropy action of $h_{s,u}$ on each component $su(2)^i$ contains the sum of three trivial
representations. Now look at the table of the irreducible 7-dimensional representation $\pi_1$
of $g_2$ (see e.g. [13, table 1]) we know that the weights of this representation are $\pm \eps_i, 0$. Here
$i = \overline{1, 3}$ and  $\eps_i-\eps_j$ $\pm \eps_i$ are the roots of $g_2$. The complexification of the
maximal compact algebra $so(4)$ in $g_2$ is the direct sum $so(4)_{\otimes \C} =< h_{\eps_1} , e_{\eps_1} , e_{-\eps _1} >_{\otimes \C} \oplus 
<h_{\eps _2 - \eps_3}, e_{\eps_2 -\eps _3}, e_{\eps_3 -\eps_2} >_{\otimes \C}$. Hence the restriction of the 
$\pi_1$ to $sl(2) =< h_{\eps_1}, e_{\eps_1}, e_{-\eps_1} >_{\otimes \C}$ is
the sume of two irreducible components of dimension 2 and the adjoint representation of
dimension 3. Thus it has no trivial component. We arrive at a contradiction. This proves
Proposition 3.2.\QED

\section{Appendix. A necessary and sufficient condition for a closed
7-manifold to admit a 3-form of $\tilde G_2$-type.}

{\bf A.1. Theorem.} {\it  Suppose that $M^7$ is a compact 7-manifold. Then $M^7$ admits a 3-form of
$\tilde G_2$-type, if and only if $M^7$ is orientable and spinnable. Equivalently the first and the second
Stiefel-Whitney classes of $M^7$ vanish.}

\medskip

{\it Proof.} If $M^7$ admits a $\tilde G_2$-structure, it admits also a $G_2$-structure, since a maximal
compact group of $\tilde G_2$ is also a  compact subgroup of the group $G_2$. Applying the Gray
criterion for the existence of a $G_2$-structure [9] we obtain the ``only if" statement. Now let
us prove the ``if " part. By a result of Dupont [4] any compact orientable 7-manifold admits
three linearly independent vector fields. Hence $M^7$ admits an $SO(4)$-structure. In particular $M^7$
admits an $SO(3) \times SO(4)$-structure. Denote by $Spin(3, 4)$ the Lie subgroup in $Spin(7,\C)$
whose Lie algebra is $so(3, 4)$. It is known (see e.g. [2 Theorem5]) that $\pi_1(Spin(3, 4)) = \Z_2$.
Denote by $(SO(3) \times  SO(4))^*$ the connected subgroup in $Spin(3, 4)$ whose Lie algebra is
$so(3)\times so(4)$. Clearly $(SO(3)\times SO(4))^*$ is a maximal compact Lie subgroup of $Spin(3, 4)$.
Taking into account the isomorphism $Spin(7)/(SO(3)\times SO(4))^* = SO(7)/(SO(3)\times SO(4))$
we conclude that $M^7$ admits an $(SO(3)\times SO(4))^*$-structure and hence a $Spin(3, 4)$-structure.
Now we shall prove that the $Spin(3, 4)$-structure on $M^7$ is reduced to a $\tilde G_2$-structure. It
is easy to see that the quotient $Spin(3, 4)/\tilde G_2$ is the pseudo-sphere $S^7(4, 4)$ in the space
$\R e_0\times  \R ^7$ of the spinor representation of $Spin(3, 4)$. This pseudo-sphere bundle over $M^7$
admits a section $e_0$ . Hence $M^7$ admits a $\tilde G_2$-structure. \QED

Using Theorem A.1 and the Gromov h-principle (see e.g. [5]) for open 3-forms of $\tilde G_2$-type we can get a lot examples of open 7-manifolds admitting a closed 3-form of $\tilde G_2$-type.
\medskip

 {\bf Acknowledgement}. I am thankful to Jiri Vanzura and Jarolim Bures for their introducing me to the field of multi-symplectic 3-forms.  The main part of this paper has been written at IHES.   This paper is partially supported by contract RITA-CT-2004-505493 and
by  grant  of  ASCR Nr IAA100190701.
 
 %\medskip

\medskip

address: Mathematical Institute  of ASCR,
Zitna 25, Praha 1, CZ-11567 Czech Republic,
  email: hvle@math.cas.cz

\end{document}